\documentclass[12pt]{article}
\usepackage[cp1251]{inputenc}
\usepackage[russian,english]{babel}
\usepackage{latexsym,amsfonts,amssymb,amsmath,amsfonts}
\usepackage{setspace}

 \newtheorem{theorem}{Theorem}[section]

 \newtheorem{remark}[theorem]{Remark}

 \numberwithin{equation}{section}

\begin{document}


\begin{center}

{\bf  Identification of the  Heat Transfer Coefficient Using an Inverse Heat Conduction Model}

\medskip
S.G. Pyatkov \footnote{ MSC Primary 35R30; Secondary 35R25;
35K57;

}

\bigskip

\parbox{10cm}{\small
 Abstract.   Inverse problems of recovering heat transfer coefficient 
 from integral measurements  are considered. The heat transfer coefficient occurs in the transmission conditions of imperfect contact type or 
  the Robin type boundary conditions.  
  It is representable as a finite part of the Fourier series with time dependent coefficients.  The additional measurements  are integrals of a solution multiplied by some weights.
 Existence and uniqueness of solutions in Sobolev classes are proven and the conditions on the data are sharp. These conditions include  smoothness and consistency conditions on the data and additional conditions on the kernels of the integral operators used in additional measurements.  The proof relies on
a priori bounds and the contraction mapping principle. The existence and uniqueness theorems are local in time.
 }

\end{center}



\section*{Introduction}
\hspace{0.7 cm} Under consideration is a parabolic equation of the form
 \begin{equation}\label{e1}
M u= u_t-Lu=f,\  Lu= \sum_{i,j=1}^{n}a_{ij}(t,x)u_{x_i x_{j}}-\sum_{i=1}^{n}a_{i}(t,x)u_{x_{i}}-a_{0}(t,x)u,
\end{equation}
where $ x\in G$ and $G\subset {\mathbb R}^{n}$ is a bounded domain with boundary
 $\Gamma$ of class  $C^{2}$ (see the definitions in
\cite[Ch.~1]{lad}), $t\in (0,T)$. Let  $Q=(0,T)\times G$,
$S=(0,T)\times \Gamma$.

This equation is a vital tool in scientific and engineering applications to assess and forecast temperature changes over time.
According to Animasaun I.L. et al. (2022) \cite{ani},
it is commonly used to model heat conduction, diffusion, and numerous dynamic thermal processes.
The problems of identification of the  heat transfer coefficients  arise in various problems of mathematical physics  (see \cite{ali,ozi,bae,tka}):
 diagnostics and identification of heat transfer in supersonic heterogeneous
flows, modeling and description of heat transfer in heat-shielding materials
and coatings,  thermal protection design  and control of heat transfer regimes, modeling of properties and thermal processes  in reusable thermal protection
of aerospace vehicles,  composite materials, ecology,  etc.

Two inverse problems are examined. In the former case, the heat transfer coefficient defined on the interface occurs in  a transmission
condition   of imperfect contact type and in the latter in the Robin boundary condition. The statements of the problems are as follows. In the former case, the domain  $G$ is divided into two open sets  $G^+$ and
    $G^-$,  $\overline{G^-} \subset G$,
    $\overline{G^+} \cup \overline{G^-} = \overline{G}, G^+ \cap G^- = \emptyset.$
    Let  $\Gamma_0 = \partial G^+ \cap \partial G^-$, $S_0=\Gamma_0 \times(0,T)$.
    The equation  \eqref{e1} is supplemented with the initial
and boundary conditions
\begin{equation}\label{e2}
B(t,x)u|_{S} =g,\ \  u|_{t=0}=u_0(x),
\end{equation}
where $Bu=\frac{\partial u}{\partial
N}+\beta u$ or $Bu=u$,  $\frac{\partial u}{\partial
N}=\sum_{i,j=1}^{n}a_{ij}(t,x)u_{x_{j}}(t,x)n_{i}$, with
$\vec{n}=(n_{1},n_{2},\ldots,n_{n})$  the outward unit normal to  $S$, and the transmission conditions
       \begin{equation}\label{m1}
    \frac{\partial u^+}{\partial N}(t,x) - \sigma(t,x)(u^+(t,x) - u^-(t,x)) = g^+(t,x), \ \ (t,x)\in S_0,
    \end{equation}
    \begin{equation}\label{m2}
    \frac{\partial u^-}{\partial N}(t,x) =\frac{\partial u^+}{\partial N}(t,x), \ \ (t,x)\in S_0,
    \end{equation}
    where $\frac{\partial u^\pm}{\partial N}(t,x_{0}) = \lim_{x\in G^\pm, \ x\rightarrow x_0 \in \Gamma_0}
    \sum_{i,j=1}^{n} a_{ij} u_{x_{i}} \nu_j $ ($\nu$ is the unit outward normal to $\partial G^-$) and
    $u^\pm(t,x_{0}) = \lim_{x\in G^\pm, \ x\rightarrow x_0 \in \Gamma_0} u(t,x)$.
The inverse problem is to determine  a solution
 $u$ to the problem  \eqref{e1}-\eqref{m2} and the  heat transfer coefficient  $\sigma=\sum_{i=1}^{m}q_{i}(t)\Phi_{i}(t,x)$,
where the functions  $q_{i}$ are unknowns and $\{\Phi_{i}(t,x)\}$ are some basis functions. It is naturally to assume that they depend only on  $x$ but  for the sake of generality we take them  depending on all variables. 

In the latter case, the inverse problem is to determine  a solution
 $u$ to the problem  \eqref{e1}-\eqref{e2} and the function  $\beta=\sum_{i=1}^{m}q_{i}(t)\Phi_{i}(t,x)$,
where the functions  $q_{i}$ are unknowns and $\{\Phi_{i}(t,x)\}$ are some basis functions. 

In both problems,  the additional integral measurements to determine  the coefficients $\sigma$  or  $\beta$   look as follows:
\begin{equation}\label{e4}
\int_{G}u(t,x)\varphi_k(x)dx=\psi_{k}(t),\ k=1,2,\ldots,m.
\end{equation}

The transmission conditions \eqref{m1}, \eqref{m2}
agree with the conventional imperfect contact condition
at the interface (see  \cite{bae}).  If  $\sigma\to \infty$ then
we come to the diffraction problem
(see \cite[Chapter 3, Sect. 13]{lad}) in which
$u^+ =    u^-$ and     $\frac{\partial u^+}{\partial N} = \frac{\partial u^-}{\partial N}$ on $S_0$.

At present, there are many publications on the numerical solution of
the problems of the type \eqref{e1}-\eqref{e4} or \eqref{e1}, \eqref{e2}, \eqref{e4}
in the various statements.
The most usable statement provides the pointwise additional measurements, in this case the condition \eqref{e4} is replaced with the conditions
$u(t,b_j)=\psi_j(t)$ $(j=1,2,\ldots,m,\ b_j\in G)$.
It is often the case when the coefficient $\sigma$
depends only on time \cite{tka,hua,zhan,moh}  or space variables \cite{abr,abr1, duda,abr12}
(see also the bibliography and the results in  \cite{art}-\cite{zhu}).
In almost all papers, the problem is reduced to some optimal control problem and
the minimization of the corresponding quadratic functional (see \cite{tka,art,zhan,dre,hua,abr,abr1}.
Let us describe some of the already addressed results.
In the case of a sole space or time variable, the heat transfer coefficient depending on the temperature
is recovered numerically with the use of pointwise measurements in \cite{tka}.
In \cite{art} the authors determine  the heat transfer coefficients
that depend in a special manner on the additional parameters from a collection of values of a solution at given points.
In \cite{lou,abr} the Monte-Carlo method is employed to restore the heat transfer coefficient
depending on two space variables. The values of a solution
on a part of the boundary serve as the overdetermination conditions.
The simultaneous recovering of a coefficient in a parabolic equation and the heat transfer coefficient is realized in \cite{hua}.
 The pointwise overdetermination conditions
are also used in \cite{dre}, \cite{zhu}. In \cite{zhu}
 under consideration is a one-dimensional inverse problem of simultaneous
recovering the heat flow on one of the lateral boundaries and the thermal contact resistance at the interface.
The article \cite{abr1} implements the numerical determination of the heat transfer coefficient from
measurements on the available part of the outer boundary of the domain.

Several existence results are known if  the pointwise ovedetermination conditions are used instead of those in \eqref{e4}.
 If the measurement points lie on the boundary of the domain and the heat transfer coefficient occurring in the boundary condition is determined then the existence and uniqueness theorems  can be found in \cite{pya1}, \cite{kos1}, \cite{kos2}.
 The same results were obtained if the  measurement points lie at the interface. The inverse problem of determination of the interface heat transfer coefficient under certain conditions is  well-posed and  the most general existence and uniqueness theorems
 can be found in  \cite{bel1,bel2}. If the  measurement points lie in $G$ then the problem becomes ill-posed.
 The conditions \eqref{e4} were used in \cite{pya11} and \cite{kozh} to determine  the heat flux on the outer boundary and existence and uniqueness theorems are proven.
   It is often the case when the integrals in \eqref{e4} are taken over the boundary of a domain
 \cite{lesn}-\cite{sil} and the heat transfer coefficient depending on time or  space variables is determined. In these articles,  the problem is reduced to
some control problem which is studied theoretically and some existence theory is presented. But these control problems are not equivalent to the initial ones.

As for the problem \eqref{e1}-\eqref{e4} of recovering the interface heat transfer coefficient $\sigma$ and the problem \eqref{e1}, \eqref{e2}, \eqref{e4} of recovering the coefficient $\beta$,  there are no theoretical results on solvability or uniqueness of solutions to this problem
in the literature except for our articles  \cite{pya00}, \cite{pya01}.
In contrast to other articles, we look for the heat transfer coefficient in the form of a finite segment of the Fourier series and
this statement allows to obtain an approximation to the heat transfer coefficient depending on all variables and the accuracy of determination depends on just a number of measurements. This article is actually a survey of the results obtained in the articles  \cite{pya00}, \cite{pya01}. 
 The conditions on the data are described which allow to state that there are existence and uniqueness theorems  in Sobolev classes for 
 solutions to the above problems.     These conditions include  smoothness and consistency conditions on the data and additional conditions on the kernels of the integral operators used in additional measurements.  The proof relies on
a priori bounds and the contraction mapping principle. The existence and uniqueness theorems are local in time.

   \section{Preliminaries}
\hspace{0.7 cm}

The Lebesgue spaces  $L_p(G;E)$  and the Sobolev spaces $W_p^s(G;E)$, $W_p^s(Q;E)$ of vector-valued functions
taking the values  in a Banach space  $E$  (see the definitions in  \cite{tri}, \cite{denk})  are used in the article.
 The Sobolev spaces are denoted by  $W_p^s(G)$, $W_{p}^{s}(Q)$, etc., whenever $E={\mathbb R}^n$.
The inclusion  $u=(u_1,u_2,\ldots, u_k)\in W_p^s(G)$ for a vector-function means that
every of the component $u_i$ of  $u$ belongs to
  $W_p^s(G)$. By a norm of a vector, we mean the sum of the norms of its coordinates.
The H\"{o}lder  spaces $C^\alpha(\overline{G})$, $C^{\alpha,\beta}(\overline{Q}), C^{\alpha,\beta}(\overline{S})$ are defined in \cite{lad} (see also \cite{tri}).
Given an interval  $J=(0,T)$, put
$W_p^{s,r}(Q)=W_p^{s}(J;L_p(G))\cap L_p(J;W_p^r(G)$ and  $W_p^{s,r}(S)=W_p^{s}(J;L_p(\Gamma))\cap
L_p(J;W_p^r(\Gamma))$.  All coefficients of $L$ are real as well as the corresponding function spaces.

To simplify the exposition, we suppose  below that $p>n+2$.
 Denote  $(u,v)=\int_{G}u(x)v(x)dx$.
  Introduce the notations
       $Q^{\tau}=(0,\tau)\times G$,    $S^{\tau}=(0,\tau)\times \Gamma$,
       $S_{0}^{\tau}=(0,\tau)\times \Gamma_{0}$, $Q^{\pm}=(0,T)\times G^{\pm}$, $Q_{\pm}^{\tau}=(0,\tau)\times G^{\pm}$.
   Let  $B_{\delta}(b)$ be a ball centered at $b$ of radius $\delta$. The symbol $\rho(X,M)$ stands for the distance between the sets $X,M\subset {\mathbb R}^n$.

\section{Identification  of the interface heat transfer coefficient}

Describe the conditions on the data ensuring solvability of the problem.
 The operator  $L$ is assumed to be elliptic, i.~e., there exists a constant  $\delta_{0}>0$ such that
\begin{equation}\label{e7}
\sum\nolimits_{i,j=1}^{n}a_{ij}(t,x)\xi_{i}\xi_{j}\geq \delta_{0}|\xi|^{2}\  \forall \xi\in {\mathbb R}^{n},\  \forall (t,x)\in Q.
\end{equation}
 The conditions on the coefficients are as follows:
 \begin{equation}\label{e8}
      a_{i}\in L_{p}(Q)\ (i\geq 0),\   a_{ij}\in
    C(Q^{\pm})\ (\ i,j=1,\ldots,n);
    \end{equation}
  the functions  $a_{ij}|_{Q^{\pm}}$ admits extensions to continuous functions of class
     $C(\overline{Q^{\pm}})$ and
    \begin{equation}\label{e9}
a_{ij}^{\pm}|_{S_{0}}\in W_p^{s_{0},2s_{0}}(S_{0}),\  a_{ij}|_{Q^{\pm}}\in C([0,T];W_p^1(G^\pm)),\  a_{ij}|_{S}\in
W_p^{s_{0},2s_{0}}(S),
    \end{equation}
where $ \ i,j=1,\ldots,n,$ $a_{ij}^\pm(t,x_{0}) = \lim_{x\in G^\pm, \ x\rightarrow x_0\in \Gamma_0} a_{ij}(t,x)$, the last inclusion in \eqref{e9} is fulfilled provided that  $Bu\neq u$ in \eqref{e2};
\begin{equation}\label{e10}
a_{ij}, a_{k} \in L_{\infty}(G;W_{p}^{s_0}(0,T))\  (k=0,1,\ldots,n, \ i,j=1,\ldots,n).\
 \end{equation}
 The main conditions on the data are of the form
\begin{equation}\label{e11}
f\in L_{p}(Q),\ \ u_{0}(x)\in W_{p}^{2-2/p}(G^{\pm}),\  g\in W_{p}^{k_0,2k_0}(S),\ g^{+}\in W_{p}^{s_0,2s_{0}}(S_{0}),
\end{equation}
where  $k_{0}=1-1/2p$ in the case of $Bu=u$ and $k_{0}=1/2-1/2p$ otherwise;
\begin{equation}\label{e12}
  \beta\in W_{p}^{s_0,2s_{0}}(S),\  g(0,x)|_{\Gamma}=B(0,x)u_{0}|_{\Gamma},\
\frac{\partial u_{0}^{+}}{\partial N}=\frac{\partial u_{0}^{-}}{\partial N},\  x\in \Gamma_{0}; \
\end{equation}
\begin{multline}\label{e202}
   \varphi_{k}|_{G^\pm}\in W_{\infty}^{1}(G^\pm),\ \Phi_{k}\in W_{p}^{s_0,2s_0}(S_0),\ \psi_{k}\in W_{p}^{s_0+1}(0,T), \\ (f,\varphi_{k})\in W_{p}^{s_0}(0,T),\ k=1,2,\ldots,m.
 \end{multline}
Assume that a pair $(u,\vec{q})$, $\vec{q}=(q_1,q_2,\ldots,q_m)$ is a solution to the problem \eqref{e1}-\eqref{e4}.
Multiply \eqref{e1} by $\varphi_i$ and integrate over $G$. Integrating by parts and using the transmission conditions, we infer
\begin{multline}\label{e21}
\psi_i'(t)+\sum_{j=1}^mq_j(t)\int_{\Gamma_0} \Phi_j(t,x)(u^+(t,x)-u^-(t,x))(\varphi^+_i(x)-\varphi^-_i(x))d\Gamma_0
- \\  \int_\Gamma \frac{\partial u}{\partial N}\varphi_i(x)\,d\Gamma 
+\int_{\Gamma_0} g^+(\varphi^+_i(x)-\varphi_i^-(x))\,d\Gamma +a(u,\varphi_i)=\int_G f(t,x)\varphi_i(x)\,dx.\\
a(u,\varphi_i)(t)=\sum_{k,l=1}^n\int_G a_{kl}u_{x_l}\varphi_{i x_k }(x)\,dG
+ \int_G (\sum_{k=1}^n a_{k}u_{x_k}+a_0u)\varphi_{i}(x)\,dG,
\end{multline}
where $\varphi_k^\pm(x_0)=\lim_{x\to x_0, x\in G^\pm} \varphi_k(x)$. Define the function  $\varphi_i^0(x)=\varphi^+_i(x)-\varphi_i^-(x)$ $(x\in \Gamma_0)$.
We would like to have that the system  \eqref{e21} is uniquely solvable relatively the vector-function $\vec{q}$, i.~e.,
$
|{\rm det\,}B(t)|\geq \delta_{0}>0 \ \ \forall t\in [0,T],
$
where $B(t)$ is the matrix with entries
$\int_\Gamma \Phi_j(t,x)(u^+(t,x)-u^-(t,x))(\varphi^+_i(x)-\varphi_i^-(x))\,d\Gamma$. Let  $B_0=B(0)$.
Taking  $t=0$, we obtain the condition
\begin{equation}\label{e22}
|{\rm det\,}B_0|\neq 0, \ b_{ij}=\int_\Gamma \Phi_j(0,x)(u^+_0(x)-u_0^-(x))(\varphi^+_i(x)-\varphi_i^-(x))\,d\Gamma.
\end{equation}
Let  $t=0$ in \eqref{e21}. We arrive at the system
\begin{multline}\label{e23}
\psi_i'(0)+\sum_{j=1}^mq_j(0)\int\limits_{\Gamma_0} \Phi_j(0,x)(u_0^+(x)-u_0^-(x))\varphi_i^0(x)\,d\Gamma_0
-  \\ \int\limits_\Gamma \frac{\partial u_0}{\partial N}\varphi_i(x)\,d\Gamma
+\int\limits_{\Gamma_0} g^+(0,x)\varphi_i^0(x)\,d\Gamma +a(u_0,\varphi_i)=(f(0,x),\varphi_i),
\end{multline}
where $\  i=1,2,\ldots,m. $
Under the condition \eqref{e22},
there exists a unique solution  $(q_1(0),\ldots,q_m(0))$ to the system \eqref{e23}.
Thus, we have determined  the function $\sigma(0,x)=\sum_{i=1}^m q_i(0)\Phi_i(0,x)$.
Taking $t=0$ at \eqref{m1}, \eqref{e4} and using the initial conditions \eqref{e2}, we come to the necessary consistency conditions
 \begin{equation}\label{e24}
\frac{\partial u_0^+}{\partial N}-\sigma(0,x)(u_0^+-u_0^-)\bigr|_{\Gamma} =g^+(0,x), \
\int_{G}u_{0}(x)\varphi_{k}(x)\,dx=\psi_{k}(0),
\end{equation}
where $k=1,\ldots,m$.
The main result of this section is the following theorem.

\begin{theorem}\label{main} Let the conditions  {\rm \eqref{e7}-\eqref{e202},   \eqref{e22}, \eqref{e24}} hold. Then, on some segment   $[0,\tau_0]$ $(\tau_0\leq T)$, there exists a unique  solution
  $(u,\vec{q})$ $(\vec{q}=(q_{1},\ldots,q_{m}))$ to the problem  {\rm \eqref{e1}--\eqref{e4}} such that
$u|_{Q^\pm}\in W_{p}^{1,2}(Q^{\tau_0}_\pm)$, $\vec{q}\in W_p^{s_0}(0,\tau_0)$.
\end{theorem}

{\bf Proof}. Outline the proof (see \cite{pya00}).  Let a pair $u\in W_{p}^{1,2}(Q^+)\cap W_{p}^{1,2}(Q^-)$, $\vec{q}\in W_{p}^{s_0}(0,T)$ be a solution to the problem  \eqref{e1}-\eqref{e4}. As before,  we can find constants
$q_{i}(0)$.
Let $\sum_{i=1}^{m}q_{i}(0)\Phi_{i}(t,x)=\sigma_{0}(t,x)$ and denote by
 $v\in W_{p}^{1,2}(Q^+)\cap W_{p}^{1,2}(Q^-)$ a solution to the auxiliary transmission  problem 
 \begin{equation}\label{e5}
    Mu = f(t,x),\ (t,x)\in Q, \ Bu|_S = g, \ u|_{t=0} =u_{0},
        \end{equation}
    \begin{equation}\label{e6}
     B^{+}u= \frac{\partial u^{+}}{\partial N}-{\sigma}_0(u^{+}-u^{-})= g^+,\
\frac{\partial u^{+}}{\partial N}=\frac{\partial u^{-}}{\partial N},\ (t,x)\in S_{0} \
    \end{equation}
          whose solvability is established with the use of Theorem 1 in \cite{bel1}.
 Make the change of variables  $u=v+w$. Inserting this function $u$ in \eqref{e1} and involving the equation  \eqref{e5}, we obtain that the function  $w\in W_{p}^{1,2}(Q^+)\cap W_{p}^{1,2}(Q^-)$ is a solution to the problem
\begin{multline}\label{o72}
w_t-Lw=0,\ \ Bw|_\Gamma=0,\ \frac{\partial w^+}{\partial N}=\frac{\partial w^-}{\partial N},\
 w|_{t=0}=0,\\
 \frac{\partial w^+}{\partial N}-\sigma_0(w^+-w^-)=(\sigma-\sigma_0)(v^++w^+-v^--w^-).
\end{multline}
The condition  \eqref{e4} is rewritten as follows:
\begin{equation}\label{o82}
\int_{G} w\varphi_{k}(x)\,dx=\psi_{k}-\int_{G}v(t,x)\varphi_{k}(x)\,dx=\tilde{\psi}_{k},\  k=1,2,\ldots,m.
\end{equation}
In view of \eqref{e202} and \eqref{e24}, $\tilde{\psi}_{k}(0)=0$ and
 $\tilde{\psi}_{k}(t)\in W_{p}^{1}(0,T)$.
  Multiply  the equation in \eqref{o72} by  $\varphi_{k}(x)$ and integrate over
 $G$. Integrating by parts yields 
\begin{multline}\label{o9}
\tilde{\psi}_i'(t)+\sum_{j=1}^m\tilde{q}_j(t)\int\limits_{\Gamma_0} \Phi_j(t,x)(w^+(t,x)-w^-(t,x)+v^+(t,x)-v^-(t,x))\varphi_i^0(x)\,d\Gamma_0-
 \\  \int\limits_\Gamma \frac{\partial w}{\partial N}\varphi_i(x)\,d\Gamma
+\int\limits_{\Gamma_0} \sigma_0(w^+(t,x)-w^-(t,x))\varphi_i^0(x)\,d\Gamma +a(w,\varphi_i)=0,
\end{multline}
where $i=1,\ldots,m$ and $\tilde{q}_i=q_i-q_i(0)$.
The equality \eqref{o9} is rewritten as
\begin{multline*}
\sum_{j=1}^m\tilde{q}_j(t)\int_{\Gamma_0} \Phi_j(t,x)(v^+-v^-)\varphi_i^0(x)d\Gamma_0
=-a(w,\varphi_i)+\int_\Gamma \frac{\partial w}{\partial N}\varphi_i(x)\,d\Gamma- \\ \tilde{\psi}_i'(t)+\int_{\Gamma_0} \sigma_0(w^+-w^-)\varphi_i^0(x)\,d\Gamma-
\sum_{j=1}^m\tilde{q}_j(t)\int_{\Gamma_0} \Phi_j(t,x)(w^+-w^-)\varphi_i^0(x)d\Gamma_0.
\end{multline*}
and, thereby, we have the operator equation
\begin{multline*}
B(t)\vec{q}=\vec{F},\  \  F_{k}=
-a(w,\varphi_i) -\tilde{\psi}_i'(t) -  \int_{\Gamma_0} \sigma_0(w^+-w^-)\varphi_i^0(x)\,d\Gamma \\ 
+\int_\Gamma \frac{\partial w}{\partial N}\varphi_i(x)\,d\Gamma +\sum_{j=1}^m\tilde{q}_j(t)\int_{\Gamma_0} \Phi_j(t,x)(w^+-w^-)\varphi_i^0(x)\,d\Gamma_0,
\end{multline*}
where $\vec{F}=(F_{1},\ldots,F_{m})^{T}, $ $\vec{q}=(\tilde{q}_{1},\ldots,\tilde{q}_{m})^{T}$ and $B(t)$ is the matrix with entries
$b_{ij}=\int_{\Gamma_0} \Phi_j(t,x)(v^+(t,x)-v^-(t,x))\varphi_i^0(x)\,d\Gamma_0$. Moreover, $B(0)=B_0$ and the matrix $B_0$ is nondegenerate.
The embedding theorems imply that $v\in C(\overline{Q})$, $\Phi_i\in C(\overline{S})$  (even more $v\in C^{1-(n+2)/2p,2-(n+2)/p}(\overline{Q})$)
and thereby there exist  parameters $\tau_0$ and $\delta_1>0$ such that
\begin{equation*}
|det\,B(t)|\geq \delta_1 \ \ \forall t\in [0,\tau_0].
\end{equation*}
For $\tau\leq \tau_0$,  we have that
\begin{equation}\label{e28}
\vec{q}=B^{-1}\vec{F} = R(\vec{q})=\vec{g}_0+R_{0}(\vec{q}),
\end{equation}
where  $\vec{g}_0=B^{-1}\vec{\Psi}$ and the $k$th coordinate  $\Psi_{k}$ of the vector  $\vec{\Psi}$ is of the form  $\Psi_{k}(t)=-\tilde{\psi}_{k}'(t)$.
This equation is used to determine  $\vec{q}$.
It is not difficult to demonstrate  that   $\int_{G}v_{t}(t,x)\varphi_{k}(x)\,dx\in W_{p}^{s_0}(0,T)$
  and, thus,  $\tilde{\psi}_{k}(t)\in W_{p}^{1+s_0}(0,T)$,  $\tilde{\psi}_{k}(0)=\tilde{\psi}_{k}'(0)=0$.
 Next, using   the conventional estimates for solutions to parabolic problems, we can show  that  
 the operator $R$ is a contraction in the ball  $B_{R_0}=\{\vec{q}\in \tilde{W}_p^{s_0}(0,\tau):\
\|\vec{q}\|_{ \tilde{W}_p^{s_0}(0,\tau)}\leq R_0\}$ and takes this ball into itself provided that the parameter $\tau$ is sufficiently small, where
$R_0=2\|\vec{g}_0\|_{\tilde{W}_p^{s_0}(0,T)}$. The contraction mapping principle implies  that the equation \eqref{e28}
is solvable locally in time (see the complete proof in \cite{pya00}).  Thus, the equation \eqref{e28} is solvable and we have determined 
 the vector $\vec{q}$. A solution $w$ in this case is a solution to  the problem \eqref{o72}.
Validate  the conditions \eqref{o82}.
Multiply the equation in  \eqref{o72} by $\varphi_{k}$ and integrate the result over  $G$.
Integrating by parts yields 
\begin{multline*}
\int\limits_{G}w_{t}\varphi_{k}\,dx+\sum_{j=1}^m\tilde{q}_j(t)\int_{\Gamma_0} \Phi_j(t,x)(w^+(t,x)-w^-(t,x)+v^+(t,x)-v^-(t,x))\varphi_i^0(x)\,d\Gamma_0
 \\ - \int_\Gamma \frac{\partial w}{\partial N}\varphi_i(x)\,d\Gamma
+\int_{\Gamma_0} \sigma_0(w^+-w^-)\varphi_i^0(x)\,d\Gamma +a(w,\varphi_i)=0.
\end{multline*}
Subtracting this equality from \eqref{o9}, we infer
$$
\int\limits_{G}w_{t}\varphi_{k}\,dx=\tilde{\psi}_{k}', \ \ k=1,\ldots,m.
$$
Integrating this equality with respect to $t$, we establish  \eqref{o82}. 
The uniqueness of solutions follows from  the estimates obtained in the proof and standard arguments.

\begin{remark} Generally speaking, the interface $\Gamma_0$ as well as the outer boundary $\Gamma$
can consist of several connectedness components.  In particular,  we can have several heat transfer coefficients occurring in different transmission conditions.
 The claim of the theorem remains valid under the same conditions. 
 \end{remark}

 \section{Identification of the heat transfer coefficients in the Robin boundary condition }
  
 The problem  \eqref{e1}, \eqref{e2}, \eqref{e4} of recovering the coefficient $\beta$ in the Robin boundary condition is considered.
Proceed with the conditions on the data of the problem. They are quite similar to those in the previous section.  
 The conditions on the coefficients are as follows:
\begin{equation}\label{e17}
a_{ij}\in C([0,T];W_p^1(G)), \  a_{ij}|_S\in W_p^{s_0,2s_0}(S)\ (s_0=1/2-1/2p),
 \end{equation}
\begin{equation}\label{e171}
a_{ij}, a_{k} \in L_{\infty}(G;W_{p}^{s_0}(0,T))\  (k=0,1,\ldots,n, \ i,j=1,\ldots,n,\  p>n+2).  
 \end{equation}
The main conditions on the data of the problem have the form  
\begin{equation}\label{e19}
f\in L_{p}(Q),\ \ u_{0}(x)\in
W_{p}^{2-2/p}(G),\ g\in W_{p}^{s_0,2s_0}(S). 
\end{equation}
Write out the additional conditions 
\begin{multline}\label{e193}
   \varphi_{k}\in W_{\infty}^{1}(G),\ \Phi_{k}\in W_{p}^{s_0,2s_0}(S),\ \psi_{k}\in W_{p}^{s_0+1}(0,T), \\ (f,\varphi_{k})\in W_{p}^{s_0}(0,T),\ k=1,2,\ldots,m.
 \end{multline}
Assume that a pair $u\in W_p^{1,2}(Q)$, $\vec{q}=(q_1,q_2,\ldots,q_m)$ is a solution to the problem \eqref{e1}, \eqref{e2}, \eqref{e4}.
Multiply \eqref{e1} by $\varphi_i$ and integrate over $G$. Integrating by parts, we infer
\begin{multline}\label{o1}
\psi_i'(t)+\sum_{j=1}^mq_j(t)\int_\Gamma \Phi_j(t,x)u(t,x)\varphi_i(x)\,d\Gamma+
\sum_{k,l=1}^n\int_G a_{kl}u_{x_l}\varphi_{x_k i}(x)\,dG \\
- \int_\Gamma g(t,x)\varphi_i(x)\,d\Gamma+ \int_G (\sum_{k=1}^n a_{k}u_{x_k}+a_0u)\varphi_{i}(x)\,dG=(f,\varphi_i).
\end{multline}
It is naturally to assume that this system  is uniquely solvable relatively the vector-function $\vec{q}$.
Thus, it is desirable to have that 
$
|{\rm det\,}B(t)|\geq \delta_{0}>0 \ \ \forall t\in [0,T],
$
where $B(t)$ is the matrix with entries 
$\int_{\Gamma}\varphi_{i}(x)\Phi_{j}(t,x)u(t,x)\,d\Gamma$.
At $t=0$ we must have 
\begin{equation}\label{o2}
|{\rm det\,}B_0|\neq 0 \ \ B_0=B(0), \ b_{ij}=\int_{\Gamma}\varphi_{i}(x)\Phi_{j}(0,x)u_0(x)\,d\Gamma.
\end{equation}
Taking $t=0$ in \eqref{o1}, we arrive at the system
\begin{multline}\label{o3}
\psi_i'(0)+\sum_{j=1}^mq_j(0)\int_\Gamma \Phi_j(0,x)u_0(x)\varphi_i(x)\,d\Gamma+
\sum_{k,l=1}^n\int_G a_{kl}(t,0)u_{0x_l}\varphi_{x_k i}(x)\,dG \\
- \int_\Gamma g(0,x)\varphi_i(x)\,d\Gamma+ \int_G (\sum_{k=1}^n a_{k}u_{0 x_k}+a_0u_0)\varphi_{i}(x)\,dG=(f(0,x),\varphi_i), 
\end{multline}
where $\  i=1,2,\ldots,m. $
Under the condition \eqref{o2}, 
there exists a unique solution $(q_1(0),\ldots,q_m(0))$ to this system. 
Thus, we have determined  the function $\beta(0,x)=\sum_{i=1}^m q_i(0)\Phi_i(0,x)$.
The consistency conditions at $t=0$ provide the equalities
 \begin{equation}\label{o4}
\frac{\partial u_0}{\partial N}+\beta(0,x)u_0\bigr|_{\Gamma} =g(0,x)\ (x\in \Gamma), \ 
\int_{G}u_{0}(x)\varphi_{k}(x)\,dx=\psi_{k}(0),\  \ k=1,\ldots,m;
\end{equation}
The main result of this section  is the following theorem. 

\begin{theorem}  Let the conditions  {\rm
  \eqref{e7}, \eqref{e17}-\eqref{e193},  \eqref{o2}, \eqref{o4}} hold. Then on some segment $[0,\tau_0]$ $(\tau_0\leq T)$ there exists  a unique  solution   $(u,\vec{q})$ $(\vec{q}=(q_{1},\ldots,q_{m}))$ to the problem  {\rm \eqref{e1}, \eqref{e2}, \eqref{e4}} such that 
$u\in W_{p}^{1,2}(Q^{\tau_0}),$ $\vec{q}\in W_{p}^{s_0}(0,\tau_0)$.
\end{theorem} 

 Let a pair $u\in W_{p}^{1,2}(Q)$, $\vec{q}\in W_{p}^{s_0}(0,T)$ is a solution to the problem  \eqref{e1}-\eqref{e4}. In view 
of \eqref{o2}, we can find constants 
$q_{i}(0)$ from the system \eqref{o3}.
Let $\sum_{i=1}^{m}q_{i}(0)\Phi_{i}(t,x)=\beta_{0}(t,x)$ and denote by 
 $v\in W_{p}^{1,2}(Q)$ a solution to the problem
\begin{equation}\label{o6}
v_t-Lv=f,\ \ \frac{\partial v}{\partial N}+\beta_0(t,x)v\bigr|_{\Gamma} =g(t,x),\  v|_{t=0}=u_{0}(x).
\end{equation}
Note that   $\vec{q}\in W_{p}^{s_0}(0,T)$ and  $\Phi_{j}\in W_{p}^{s_0,2s_0}(S)$ then 
 $q_{i}(t)\Phi_{i}(t,x)\in W_{p}^{s_0,2s_0}(S)$,  and  $g \in W_{p}^{s_0,2s_0}(S)$ as well.
Make the change of variables  $u=v+w$. The function  $w\in W_{p}^{1,2}(Q)$ is a solution to the problem 
\begin{equation}\label{o7}
w_t-Lw=0,\ \ \frac{\partial w}{\partial N}+\beta_0(t,x)w\bigr|_{\Gamma} =(\beta_0-\beta)(v+w),\  \omega|_{t=0}=0.
\end{equation}
The condition  \eqref{e4} is rewritten as follows: 
\begin{equation}\label{o8}
\int_{G} w\varphi_{k}(x)\,dx=\psi_{k}-\int_{G}v(t,x)\varphi_{k}(x)\,dx=\tilde{\psi}_{k},\  k=1,2,\ldots,m.
\end{equation}
In view of  \eqref{o4}, $\tilde{\psi}_{k}(0)=0$ and at least
 $\tilde{\psi}_{k}(t)\in W_{p}^{1}(0,T)$. 
 It is easy to demonstrate  that   $\tilde{\psi}_{k}(t)\in W_{p}^{1+s_0}(0,T)$
  and, thus, $\int_{G}v_{t}(t,x)\varphi_{k}(x)\,dx\in W_{p}^{s_0}(0,T)$.
  Multiply  the equation in \eqref{o7} by  $\varphi_{k}(x)$ and integrate over 
 $G$. We obtain that 
$(w_{t},\varphi_{k})=(Lw, \varphi_{k}).$ 
Integrating by parts, we infer 
\begin{equation}\label{la}
\tilde{\psi}_{k}'(t)+a(w,\varphi_{k})+ \int\limits_{\Gamma}\beta_0 w \varphi_{k}\,d\Gamma+  \sum_{i=1}^{m}\tilde{q}_{i}(t)\int\limits_{\Gamma}(v+w)\Phi_{i}\varphi_{k}\,d\Gamma=0,\
\end{equation}
where $i,k=1,\ldots,m$  and
$
a(w,\varphi_{k})=\int_{G}\sum_{i,j=1}^{n}a_{ij}\omega_{x_{j}}\varphi_{kx_{i}}+
(\sum_{i=1}^{n}a_{i}\omega_{x_{i}}+a_{0}\omega)\varphi_{k}\,dx.  
$
The last equality is rewritten as 
\begin{multline}\label{e26}
\sum_{i=1}^{m}\tilde{q}_{i}(t)\int\limits_{\Gamma}\Phi_{i}\varphi_{k}v(t,x)\,d\Gamma=
-\sum_{i=1}^{m}\tilde{q}_{i}(t)\int\limits_{\Gamma}\Phi_{i}\varphi_{k}w\,d\Gamma - \\
\tilde{\psi}_{k}'(t)-a(\omega,\varphi_{k})-\int\limits_{\Gamma}\beta_0 w \varphi_{k}\,d\Gamma
\end{multline}
and thereby 
\begin{equation}\label{e27}
B(t)\vec{q}=\vec{F},\  \vec{F}=(F_{1},\ldots,F_{m})^{T},\ \vec{q}=(\tilde{q}_{1},\ldots,\tilde{q}_{m})^{T}, 
\end{equation}
where $  F_{k}=-\sum_{i=1}^{m}\tilde{q}_{i}(t)\int_{\Gamma}\Phi_{i}\varphi_{k}w\,d\Gamma -
\tilde{\psi}_{k}'(t) -a(\omega,\varphi_{k})-\int_{\Gamma}\beta_0 w \varphi_{k}\,d\Gamma $ and $B(t)$ is the matrix with entries 
$b_{ij}=\int_{\Gamma}\Phi_{j}\varphi_{i}v(t,x)\,d\Gamma$. Note that $B(0)=B_0$ and the matrix $B_0$ is nondegenerate.
The embedding theorems imply that $v\in C(\overline{Q})$, $\Phi_i\in C(\overline{S})$  (even more $v\in C^{1-(n+2)/2p,2-(n+2)/p}(\overline{Q})$)
and thereby there exist  parameters $\tau_0$ and $\delta_1>0$ such that 
\begin{equation}\label{e271}
|det\,B(t)|\geq \delta_1 \ \ \forall t\in [0,\tau_0].
\end{equation}
The function  $w$  in  \eqref{e27} is a solution   to the problem 
 \eqref{o7}. 
For $\tau\leq \tau_0$,  we have that 
\begin{equation}\label{e28}
\vec{q}=B^{-1}\vec{F} = R(\vec{q})=\vec{g}_0+R_{0}(\vec{q}),
\end{equation}
where  $\vec{g}_0=B^{-1}\vec{\Psi}$ and the $k$-th coordinate  $\Psi_{k}$ of the vector  $\vec{\Psi}$ is of the form  $\Psi_{k}(t)=-\tilde{\psi}_{k}'(t)$.
We use this equation to determine  $\vec{q}$.
Next, using the known estimates for solutions to parabolic problems, we demonstrate that the operator $R$ is a contraction in the ball  $B_{R_0}=\{\vec{q}\in \tilde{W}_p^{s_0}(0,\tau):\ 
\|\vec{q}\|_{ \tilde{W}_p^{s_0}(0,\tau)}\leq R_0\}$ and takes it into itself if the parameter $\tau$ is sufficiently small, where 
$R_0=2\|\vec{g}_0\|_{\tilde{W}_p^{s_0}(0,T)}$. The contraction mapping principle implies the solvability of the equation \eqref{e28}. 
We have determined the vector-function $\vec{q}$ on some time segment.    A solution $w$ in this case is a solution to  the problem \eqref{o7}.
Show that the conditions \eqref{o8} hold for a solution to the problem \eqref{o7}.
Multiply the equation in  \eqref{o7} by $\varphi_{k}$ and integrate the result over  $G$.
Integrating by parts, we infer
\begin{equation}\label{50}
\int\limits_{G}w_{t}\varphi_{k}\,dx+
=-a(w,\varphi_{k})-\int\limits_{\Gamma}\beta_0 w \varphi_{k}\,d\Gamma- \sum_{i=1}^{m}\tilde{q}_{i}(t)\int\limits_{\Gamma}(v+w)\Phi_{i}\varphi_{k}\,d\Gamma,\
\end{equation}
Subtracting this equality from \eqref{la}, we conclude that
$$
\int\limits_{G}w_{t}\varphi_{k}\,dx=\tilde{\psi}_{k}', \ \ k=1,\ldots,m,
$$
Integrating this equality with respect to $t$, we establish the equality \eqref{o8}.
The uniqueness of solutions follows from    the estimates obtained in the proof.

\section{Discussion}

We consider inverse problems of recovering the heat transfer coefficient  from integral measurements.
 These problems arise in some  practical applications,
 but there are no theoretical results concerning  the existence and uniqueness questions.
 The results  can be used in developing new numerical algorithms
and provide new conditions of existence and uniqueness of solutions to these problems.
We consider a model case, but
it is clear what changes should be made in the general case for validating similar results.
The main conditions on the data are conventional.
  The proof relies on  a priori bounds and  the contraction mapping principle.
  We think that the results will allow to establish some global existence and uniqueness results based on  
  the maximum principle  and additional conditions on the data. Similar results are valid in the H\"{o}lder spaces.
But we think that from the viewpoint of applications it is better  to deal with the Sobolev spaces.


\section{Conclusions}

The existence and uniqueness theorems in inverse problems of recovering the heat transfer coefficient  from the integral measurements
are proven locally  in time. The heat transfer coefficient occurs in the transmission conditions of imperfect contact type. It is sought in the form of a finite segment of the Fourier series with coefficients depending on time.  The proof relies on  a priori bounds and  fixed point theorem.
The conditions on the data ensuring existence and uniqueness of solutions in Sobolev
classes  are  sharp. They are smoothness and consistency conditions on the
data and additional conditions on the kernels of the integral operators used in additional measurements.

\vspace{5mm}
{ \bf Acknowledgement.}
{\it The research was carried out within the state assignment of Ministry of Science and Higher Education of the Russian Federation (theme No. FENG-2023-0004,
"Analytical and numerical study of inverse problems on recovering parameters of atmosphere or water pollution sources and (or) parameters of media")
 }

\vspace{5mm} \noindent {Sergey  Pyatkov, Doctor of Sci.,  Professor, Engineering School of  Digital Technologies, Yugra State University,  Khanty-Mansiysk, Russia, E-mail: pyatkovsg@gmail.com} 

\end{document}